\documentclass[onecolumn]{autart}
\usepackage{}    % Enable this line and disable the
\usepackage{tikz} % 作图宏包
\usetikzlibrary{shapes,arrows,chains,backgrounds,fit,decorations.pathreplacing}
\usepackage[europeanresistors,americaninductors]{circuitikz}
\usetikzlibrary{shapes,arrows}
\usepackage{amsmath,bm,times}
\usepackage{mathrsfs}
\usepackage{bbm}    % Enable this line and disable the
\usepackage{lineno,hyperref}
\usepackage{graphicx}           % Include figure files
\usepackage{dcolumn}            % Align table columns on decimal point
\usepackage{bm}                 % bold math

\usepackage{graphics,subfigure,epsfig,wrapfig}
\usepackage{subfigure}
\graphicspath{{eps/}}
\usepackage{epstopdf}
\usepackage{epsfig}             % for postscript graphics files
\usepackage{times}              % assumes new font selection scheme installed
\usepackage{amsmath}
\usepackage{amssymb}
\usepackage{amsfonts}
\usepackage{fancyhdr}
\usepackage{color}
\usepackage{subfigure}
\usepackage{enumerate}
\usepackage{cite}
\usepackage{wrapfig}
\usepackage{url}
\allowdisplaybreaks[4] % [1] page breaks are allowed, but avoided if possible

\newtheorem{theorem}{Theorem} %
\newtheorem{lemma}{Lemma}

\newtheorem{remark}{Remark}

\newcommand{\diff}[1]{\operatorname{d}\!{#1}} % \diff{x}
\begin{document}

%\begin{frontmatter}
\begin{frontmatter}

\title{Input-to-state stabilization of $1$-D parabolic equations with Dirichlet boundary disturbances under boundary fixed-time control} % Title, preferably not more

 \thanks[footnoteinfo]{Corresponding author}
%\thanks{This work is supported in part by the NSFC under grant 11901482 and in part by the NSERC under grant RGPIN-2018-04571.}

\author{Jun Zheng$^{1,2,\star}$}\ead{zhengjun2014@aliyun.com},
\author{Guchuan Zhu$^{2}$}\ead{guchuan.zhu@polymtl.ca}

\address{$^{1}${School of Mathematics}, Southwest Jiaotong University, Chengdu, Sichuan 611756, P.~R.~of~China\\
                $^{2}$Department of Electrical Engineering, Polytechnique
                Montr\'{e}al,  P.O. Box 6079, Station Centre-Ville,
                Montreal, QC, Canada H3T 1J4}  % Please supply

\begin{abstract}
This paper addresses the problem of stabilization of $1$-D parabolic equations with  destabilizing terms and    Dirichlet boundary disturbances. By using the method of  backstepping and the technique of splitting, a  boundary feedback controller is designed to ensure  the input-to-state stability (ISS) of the closed-loop system with   Dirichlet boundary disturbances, while preserving fixed-time stability (FTS) of the corresponding disturbance-free system, for which the {fixed  time} is either determined by the Riemann zeta function  or freely prescribed. To overcome the difficulty brought by  Dirichlet boundary disturbances,  the ISS and  FTS  properties of the involved systems are  assessed by applying the generalized Lyapunov method. Numerical simulations are conducted to illustrate the effectiveness of the proposed scheme of control design.
%
%The aim of this paper is to design a boundary control to ensure  the input-to-state stability (ISS)
% for  a class of $1$-D parabolic equations with  destabilizing terms and    Dirichlet boundary disturbances, while preserving fixed-time stability for the corresponding disturbance-free systems, for which the fixed-time is either determined by the Riemann zeta function  or freely prescribed.
% The controller is designed by means of  backstepping combined with the technique of splitting, and the ISS with respect to the   Dirichlet boundary disturbances is assessed by applying the generalized Lyapunov method.
\end{abstract}

\vspace*{-12pt}
                              % keyword list or with the
                                           % help of the Automatica
      %\footnotetext{This work is supported in part by the NSFC under grant 11901482 and in part by the NSERC under grant RGPIN-2018-04571.}                          % keyword wizard
\end{frontmatter}
\endNoHyper
%\linenumbers

%\tableofcontents %
%%%%%%%%%%%%%%%%%%%%%%%%%%%%%%%%%%%%%%%%%%%%%%%%%%%%%%%%%%%%%%%%%%%%%%%%%%%%%%%%%%%%%%%%%%%%%%%%%

%%%%%%%%%%%%%%%%%%%%%%%%%%%%%%%%%%%%%%%%%%%%%%%%%%%%%%%%%%%%%%%%%%%%%%%%%%
\section{Introduction}\label{Sec: Introduction}
  It is well-known that the  backstepping method (see \cite{krstic2008boundary}) is  effective   for designing   boundary   controllers to stabilize partial differential equations (PDEs) with destabilizing terms; see,   e.g., {\cite{balogh2002,boskovic2001boundary,
Feng:2022,Karafyllis:2021,liu2003boundary,meurer2009tracking,
Smyshlyaev:2004,smyshlyaev2005control} for PDEs with    different types of coefficients, and    \cite{baccoli2015boundary,deutscher2018backstepping,Ghaderi:2019,Kerschbaum:2020,Kerschbaum:2022,ge2022event,Orlov:2017}} for coupled PDEs, respectively.
In addition, under the framework of  input-to-state stability (ISS)  theory of PDEs, in \cite{Karafyllis:2016a,Mironchenko:2019,Zheng:2019},   the method of backstepping was used to design
      boundary controllers to ensure  the ISS of   parabolic PDEs with destabilizing terms when   boundary   disturbances are involved.

%the input-to-state stability (ISS) theory has proven to be suitable for describing the impact of external disturbances on  stability of the systems; see, e.g., \cite{Jayawardhana:2008,Karafyllis:2018iss,mironchenko2020input}. Recently, within the framework of ISS theory,
%
%
%
%using backstepping to design controllers to ensure the ISS of PDEs with destabilizing terms  in closed loop has gradually attracted widespread attention; see, e.g., \cite{Karafyllis:2016a,Mironchenko:2019,Zheng:2019}.

 It is worth mentioning that a notable feature   in the problem of input-to-state stabilization is that the disturbance-free system in closed loop   is asymptotically stable under the designed controller. However, many control systems exhibit a property stronger than the asymptotic stability, that is the finite-time stability, which requires  the state of the  disturbance-free system to reach equilibrium within a finite time  {\cite{Golestani:2021,Haimo:1986,Holloway:2019,Li:2019,Li:2017,Meurer:2011,Smyshlyaev:2004,Song:2017,Zarchan:2012}}.

   In the past few years, the problem of finite-time stabilization  of PDEs has attracted much attention and  significant results have been obtained {\cite{bao:2022,Coron:2017,Espitia:2019,Han:2024,Ouzahra:2021,Perrollaz:2014,Polyakov:2017,Polyakov:2017b, Wei:2022,Wei:2023b,Steeves:2020,Steeves:2019,Xiao:2024}}. Notably,  by  using backstepping,   the authors of \cite{Coron:2017} designed boundary controllers to  {achieve finite-time stability of} parabolic PDEs in {the} absence of external disturbances, while the authors of \cite{Han:2024}  designed   state-dependent switching control laws to ensure the so-called dual properties of the   heat equation in closed loop, namely,   the  heat equation is finite-time stable     in {the} absence of external disturbances, while preserving the   ISS in {the} presence of  distributed in-domain disturbances. Nevertheless, for parabolic PDEs with destabilizing terms and  boundary disturbances,  the  so-called dual properties of the closed-loop system have not been studied under the backstepping control.

%considered a stabilization problem of the heat equation  with distributed in-domain disturbances. By using  backstepping, the authors successfully designed a suitable boundary feedback controller, which not only ensures that the closed-loop system is ISS, but also ensures that the disturbance-free closed-loop system is finite time stable.

  In this paper, we  study the stabilization problem of   $1$-D linear parabolic equations  involving not only distributed in-domain  but also Dirichlet boundary disturbances. Moreover, unlike \cite{Han:2024}, the open-loop system considered in this paper is allowed to be unstable.
 Based on the technique of splitting and  the method applied in \cite{Coron:2017},
  a boundary controller is designed by means of backstepping to  simultaneously ensure the ISS of the closed-loop system  and the finite-time stability  of the disturbance-free closed-loop system. In particular, the {finite time}   {is independent of initial data}.

It should be noted that, compared to    \cite{Han:2024}, the control design scheme proposed in this paper  is much easier while Dirichlet boundary disturbances present more challenges in the stability analysis   due to the fact that the classical Lyapunov method is not suitable  in the appearance of Dirichlet boundary disturbances. To avoid using state-dependent switching control laws presented in \cite{Han:2024}, we apply the technique of splitting in the control design via backstepping, while, to overcome the difficulty in the stability analysis,  we handle the terms of  Dirichlet boundary disturbances by using the generalized Lyapunov method, which was proposed in \cite[Remark 4.1]{Zheng:2024}   recently.  {More specifically, we split first the considered system    into two sub-systems, among which   one  is a disturbed sub-system   with zero initial datum while another one is a disturbance-free sub-system with non-zero initial datum. Then,
    for the disturbed sub-system,  by using  backstepping, we design a  boundary controller via a time-invariant kernel function.   By using the   generalized Lyapunov method {proposed in \cite{Zheng:2024}}, we prove that the boundedness of {the} solution  to this subsystem is determined  only by the external disturbances. For
   the disturbance-free sub-system,  {following \cite{Coron:2017}}, we design fixed-time  boundary controllers   over different time intervals, and derive sufficient conditions for time segments   to ensure the  fixed-time stability of the sub-system. In particular, the {fixed time} is  either chosen based on the Riemann zeta function or freely prescribed.  Finally,  by using  controllers of the two sub-systems,   we design a boundary controller for the original system with distributed in-domain and Dirichlet boundary disturbances to  simultaneously ensure    the  ISS of the closed-loop system  and the fixed-time stability of the  corresponding disturbance-free system.}
   
  {The main} contribution of this paper lies in
     the  design of  a boundary controller, which ensures that the parabolic equation in closed loop   is ISS with respect to (w.r.t.) Dirichlet boundary disturbances and fixed-time stable in the absence of external disturbances.

   In the rest of the paper, {we present first some basic} notations. In Section~\ref{Sec. II}, we introduce the problem formulation and provide  some preliminaries on kernel functions. In Section~\ref{Sec. IV}, we design boundary controllers and state   main results in two cases, where the finite time is determined by the Riemann zeta function and freely prescribed, respectively. In Section~\ref{Sec. V}, we assess   stability for the closed-loop system under the designed boundary control. Numerical simulations are provided in  Section~\ref{Sec. VI} to illustrate the effectiveness of the proposed scheme.
Some concluding remarks are given in Section~\ref{Sec. VII}.

\paragraph*{Notation.}
Let $\mathbb{N}$ and $\mathbb{N}_+$ be the sets of non-negative and positive integers, respectively. Let $\mathbb{R}:= (-\infty,+\infty) $,  $
 \mathbb{R}_{\geq0}:= [0,+\infty)$, $
 \mathbb{R}_{>0}:= (0,+\infty)$, and $
 \mathbb{R}_{\leq 0}:= (-\infty,0]$.  

 For $\tau\in\mathbb{R}_{>0} $, let $Q_\tau := (0,1) \times  (0,\tau)$. Let $Q_{\infty} := (0,1) \times  (0,+\infty)$ and $\overline{Q}_{\infty} := [0,1] \times  [0,+\infty)$. For   $v: D\rightarrow \mathbb{R}$ with $ D\subset \overline{Q}_{\infty}$,  the notation $v[t]$ (or $v[y]$)  denotes the profile at certain $t \in{\mathbb{R}_{\geq 0}}$, i.e., $ v[t] (y)=v(y,t)$.
  
  Throughout this paper, all the functional spaces,  as well as the  norms in normed linear spaces, are defined in the standard way as in, e.g., \cite{Evans:2010}.

For classes  of comparison functions, let
\begin{align*}
  \mathcal{P}:=&\{\gamma \in C(\mathbb{R}_{\geq 0};\mathbb{R}_{\geq 0})|\gamma(0)=0, \gamma(s)>0,\forall s\in \mathbb{R}_{>0}\},\\
 \mathcal {K}:=&\{\gamma \in \mathcal{P}|\gamma\ \text{is\ strictly\ increasing}\},\\
% \mathcal{K}_{\infty}:=&\{\gamma \in \mathcal{K}| \lim\limits_{s\to+\infty}\gamma(s)=+\infty\},\\
 \mathcal {L}:=&\{\gamma \in C(\mathbb{R}_{\geq 0};\mathbb{R}_{\geq 0})|\gamma\ \text{is\ strictly\ decreasing\ with} \ \lim\limits_{s\to+\infty}\gamma(s)=0\},\\
 \mathcal {K}\mathcal {L}:=& \{\beta\in  C(\mathbb{R}_{\geq 0}\times \mathbb{R}_{\geq 0};\mathbb{R}_{\geq 0})|\beta(\cdot,t) \in \mathcal {K}, \forall t \in \mathbb{R}_{\geq 0};  \ \beta(s,\cdot)\in \mathcal {L}, \forall s \in {{\mathbb{R}_{>0}}}\}.
 \end{align*}

\section{Problem setting and preliminaries}\label{Sec. II}
In this section, we present the problem formulation and provide some preliminaries.

\subsection{Problem setting}
 For certain fixed time $T_0>0$,   we consider the problem of input-to-state stabilization for the following $1$-D parabolic equation:
\begin{subequations}\label{original system}
\begin{align}
	 u_{t}(x,t)=& (a(x)u_{x}(x,t))_x+ c(x)u(x,t)  +f(x,t) , (x,t)\in Q_{T_0},\label{1a}\\
	 u(0,t)=&{d_0(t)}, t\in (0,T_0),\label{1b}\\
	 u (1,t)=&U(t)+d_1(t)   ,t\in (0,T_0), \label{1c} \\
	  u(x,0)=&u_0(x),x\in(0,1),
\end{align}
\end{subequations}
  where   $u_0$ is the initial {datum,}  $f$ represents {the} in-domain   disturbance,  $d_0,d_1$ represent  Dirichlet boundary disturbances, and    $U(t) $ is the control input, which will be designed later.

 Throughout this paper, we assume that $a\in H^2(0,1)$ and there exists a  constant $\Lambda\geq 1$ such that
 \begin{align*}
  \frac{1}{\Lambda}\leq a(x)\leq \Lambda,\forall x\in (0,1).
 \end{align*}
 We assume that $c\in H^1(0,1)$,  $f\in C( (0,1)\times\mathbb{R}_{\geq 0}  )$,  $d_0,d_1\in C(\mathbb{R}_{\geq 0}  )$, and  the initial datum $u_0\in L^2(0,1)$.

The   goal of this paper is to  design a  boundary control law $U(t)$ by using backstepping  such that the following so-called dual properties, which was proposed in \cite{Han:2024}, hold true {simultaneously, i.e.,}
  \begin{enumerate}[(i)]
    \item  the disturbed  system~\eqref{original system}, i.e., system~\eqref{original system} with $f^2+d_0^2+d_1^2\not\equiv0$, is input-to-state stable in the spatial $L^2$-norm w.r.t.   in-domain and boundary disturbances, namely,   there exist  functions $\beta\in \mathcal{K}\mathcal{L}$ and $\gamma,\gamma_0,\gamma_1\in \mathcal{K}$  such that
\begin{align*}
 \|u[t]\|_{L^2(0,1)}
 \leq& \beta\left( \|u_0\|_{L^2(0,1)} ,t \right) +  \gamma\left(  \|f \|_{L^\infty({Q_{t}})} \right)  +\gamma_0\left(  \|d_0 \|_{{L^\infty(0,t)}} \right) +\gamma_1\left(  \|d_1 \|_{{L^\infty(0,t)}} \right),\forall t\in (0,T_0);
\end{align*}
%holds true for all $u_0\in L^2(0,1)$;
  \item    the disturbance-free system~\eqref{original system}, i.e., system~\eqref{original system} with $f\equiv d_0\equiv d_1\equiv 0$,
      %\begin{subequations}\label{original system-2}
%\begin{align}
%	 u_{t}(x,t)=& (a(x)u_{x}(x,t))_x+ c(x)u(x,t)  ,\notag\\
%&(x,t)\in [0,1]\times (0,T),\label{2a}\\
%	 u(0,t)=&0, t\in (0,T),\label{2b}\\
%	 u (1,t)=&U(t)   ,t\in (0,T), \label{2c}\\
%	 u(x,0)=&u_0(x),x\in(0,1),
%\end{align}
%\end{subequations}
is  fixed-time stable in the spatial $L^2$-norm in time $T_0>0$, i.e., $\|u[t]\|_{L^2(0,1)}\to 0$ as $t\to T_0^-$,
 where  the  time $T_0$
 %is called fixed-time, and
      is determined by the   Riemann zeta function
  or   freely prescribed.

%where $u$ is the solution to the system~\eqref{original system} with initial data $u_0$.
\end{enumerate}

%Note that the concepts of \emph{fixed-time stability} and \emph{fixed-time} used in this paper are slight  modifications of the one given   in the existing literature; see, e.g., \cite{Han:2024}, where \emph{fixed-time stability}  is defined to be dependent of \emph{initial data} and .
%

%\begin{definition}System \eqref{original system} is said to be fixed-time input-to-state practically stable (FTISpS) if there exist $\beta\in \mathcal{G}\mathcal{K}\mathcal{L},\gamma\in \mathcal{K}$ and $\epsilon>0$ such that for all $u_0\in L^2(0,1)$ and $t\geq 0 $ it holds that
%\begin{align*}
%\|u[t]\| \leq& \beta\left( \|u_0\| ,t \right) \notag\\
%&+  \gamma\left( \sup_{s\in (0,t)}\|f(\cdot,s)\| \right)+\epsilon,\forall t\geq 0,
%\end{align*}
%\end{definition}
\subsection{Preliminaries} \label{Sec. III}
Throughout this paper, let $\lambda_0>0  $ be a fixed  constant, and define $D:=\{(x,y)\in [0,1]\times [0,1]|y\leq x\}$ and \begin{align*}
  \frac{\diff{}}{\diff{x}}g (x,x) :=g_x(x,x)+g_y(x,x) ,\forall g\in H^1(D).
 \end{align*}

We present   some   results on the existence and \emph{a priori} estimates of kernel functions, which will be used in the control design.

\begin{lemma}[{\cite[Corollary 1]{Coron:2017}}]\label{estimate-k}For every $\lambda \geq \lambda_0$, the kernel equation
\begin{subequations}\label{kernel-k}
 \begin{align}
   (a(x)k_{x}(x,y))_x-(a(y)k_{y}(x,y))_y 
 = & (\lambda+c(y))k(x,y) , \forall (x,y)\in D,\\
  2a(x)\frac{\diff{}}{\diff{x}}k(x,x)+a_x(x)k(x,x) 
  =   & -(\lambda+c(x)),\forall x\in [0,1],\\
 k(x,0)=&   0,\forall x\in [0,1],
 \end{align}
 \end{subequations}
 admits a unique solution $k\in H^1(D)$ having the estimate
\begin{align*}
  \|k\|_{H^1(D)}\leq e^{C\sqrt{\lambda}}
 \end{align*}
 for  some positive constant $C$ independent of $\lambda\in [\lambda_0,+\infty)$.
\end{lemma}
\begin{lemma}[{\cite[Corollary 2]{Coron:2017}}]\label{estimate-l}  For every $\lambda \geq \lambda_0$, the kernel equation
   \begin{subequations}\label{kernel-l}
 \begin{align}
   (a(x)l_{x}(x,y))_x-(a(y)l_{y}(x,y))_y 
   =& -(\lambda+c(y))l(x,y) , \forall (x,y)\in D,\\
  2a(x)\frac{\diff{}}{\diff{x}}l(x,x)+a_x(x)l(x,x) 
 =   &-( \lambda+c(x)),\forall x\in [0,1],\\
 l(x,0)= &0,\forall x\in [0,1],
 \end{align}
 \end{subequations}
 admits a unique solution $l\in H^1(D)$ having the estimate
 \begin{align*}%\label{estimate-l}
 \|l\|_{H^1(D)}\leq C\lambda^2
 \end{align*}
 for some positive constant $C$ independent of $\lambda\in [\lambda_0,+\infty)$.
\end{lemma}
\begin{lemma}[{\cite[Lemma 4]{Coron:2017}}]\label{k-l}    Let $k \in H^1(D)$ and $l \in H^1(D)$  be the unique solution to \eqref{kernel-k} and \eqref{kernel-l} with $\lambda \geq \lambda_0$, respectively. Let $h\in L^2(0,1)$, and for  almost  every  $x\in[0,1]$ define
\begin{align*}
\tilde{h}(x):=h(x)-\int_0^xk(x,y)h(y)\text{d}y.
\end{align*}
Then, it holds for  almost  every  $x\in[0,1]$ that
\begin{align*}
h(x):=\tilde{h}(x)+\int_0^xl(x,y)\tilde{h}(y)\text{d}y.
\end{align*}
\end{lemma}
\begin{remark}\label{Remark1}
As indicated in \cite[Chapter 4]{krstic2008boundary}, for general functions $a$ and $c$ depending on the spatial variable, there is no a closed-form solution $k$ (or $l$) to the kernel equation~\eqref{kernel-k} (or \eqref{kernel-l}). However, for  positive constants $a $ and $c$, it can be proved that (see  \cite[Chapter 4]{krstic2008boundary})
 \begin{align*}
 k(x,y)=&-\frac{\lambda+c}{a}y\frac{I_1\left(\sqrt{\frac{\lambda+c}{a}(x^2-y^2)}\right)}{\sqrt{\frac{\lambda+c}{a}(x^2-y^2)}},\\
 l(x,y)=&-\frac{\lambda+c}{a}y\frac{J_1\left(\sqrt{\frac{\lambda+c}{a}(x^2-y^2)}\right)}{\sqrt{\frac{\lambda+c}{a}(x^2-y^2)}},
 \end{align*}
 where $I_1$ is  the first-order modified Bessel
function   given by
 \begin{align*}
I_1(x):=\sum_{m=0}^{\infty}\frac{\left(\frac{x}{2}\right)^{2m+1}}{m!(m+1)!},
 \end{align*}
 and $J_1$ is  the first-order  Bessel
function   given by
 \begin{align*}
J_1(x):=\sum_{m=0}^{\infty}\frac{(-1)^m\left(\frac{x}{2}\right)^{2m+1}}{m!(m+1)!}.
 \end{align*}
\end{remark}

\section{Control design}\label{Sec. IV}
In this section, we design  boundary controllers   in two cases:

\textbf{Case I}\quad the fixed time $T_0$ is determined by the Riemann zeta function;

\textbf{Case II}\quad
  the fixed time  $T_0$ is freely prescribed.

  The main method  adopted in the control design  is backstepping applied in \cite{Coron:2017} and the technique of splitting, for which we consider $u$ as $u:=v+w$ with $v$ and $w$ being the solutions of  a disturbance-free system and a disturbed system under different boundary controls, respectively.

 \subsection{A lemma on fixed-time stability of disturbance-free system}
Throughout this paper, let $\{t_n\}_{n\in \mathbb{N}}$ be a strictly increasing sequence of real numbers such that $t_0=0$.
 %and $t_n\to T$ as $n\to \infty$.
 For   given constants  $\lambda_n\geq \lambda_0$, let $k_n\in H^1(D)$ be  the unique solution to the  kernel equation corresponding to $n\in \mathbb{N}$:
 %\begin{subequations} \label{kernel-kn}
 \begin{align*}
  (a(x)k_{n,x}(x,y))_x-(a(y)k_{n,y}(x,y))_y 
=&(\lambda_n+c(y))k_n(x,y) , \forall (x,y)\in D,\\
  2a(x)\frac{\diff{}}{\diff{x}}k_n(x,x)+a_x(x)k_n(x,x) 
=& -(\lambda_n+c(x)),\forall x\in [0,1],\\
 k_n(x,0)= &0,\forall x\in [0,1].
 \end{align*}
% \end{subequations}

Motivated by \cite{Coron:2017}, we consider  the disturbance-free system:
 \begin{subequations}\label{sub-OS-1}
\begin{align}
	 v_{t}(x,t)=& (a(x)v_{x}(x,t))_x+ c(x)v(x,t), 
  (x,t)\in [0,1]\times (t_n,t_{n+1}),\label{sub-OS-1a}\\
	v(0,t)=&0, t\in (t_n,t_{n+1}),\label{sub-OS-1b}\\
	 v (1,t)=&V(t)   ,t\in (t_n,t_{n+1}), \label{sub-OS-1c}
	 % v(x,t_n):=&v_{0n}(x),x\in(0,1),
\end{align}
\end{subequations}
with boundary control law
\begin{align}\label{control law-V}
 V(t):= \int_0^1k_n (1,y)v(y,t)\diff{y},\forall t\in (t_n,t_{n+1}),
   \end{align}
   and initial condition
  \begin{align*}v(x,t_0):=u_0(x) ,\forall x\in(0,1).\end{align*}

   Set $s_0:=0$ and
\begin{align}\label{sn}
s_n:=& {\sum\limits_{j=0}^{n-1}} \lambda_j(t_{j+1}-t_j), n\in \mathbb{N}_+.
     \end{align}

    In \cite{Coron:2017}, it was proved the following stability result, which indicates rapid convergence of solutions    to disturbance-free systems,   as well as their fixed-time stability when $n\rightarrow +\infty$.
\begin{lemma}[{\cite[Proposition~1]{Coron:2017}}]\label{lem2} For system~\eqref{sub-OS-1} under the control  law \eqref{control law-V}, there exists a positive constant $\gamma_0$ depending only on $a$ and $c$, such that if, for large $n$,
 \begin{align*}
 (t_{n+1}-t_{n}) \lambda_{n}\geq \gamma_0\sqrt{\lambda_{n+1}} ,
   \end{align*}
   then, for $t\in (t_n,t_{n+1})$,
   \begin{subequations}\label{estimate-v}
    \begin{align}
 \|v[t]\|_{L^2(0,1)}\leq& Ce^{-\frac{s_{n-1}}{4}+C(n-1)} \|u_0\|_{L^2(0,1)},\\
  |V(t) | \leq& Ce^{-\frac{s_{n-1}}{4}+C(n-1)+C\sqrt{\lambda_n}} \|u_0\|_{L^2(0,1)},
   \end{align}
   \end{subequations}
where $C$ is a positive constant independent of $n$ and $u_0$.

 In particular, if, in addition,
%\begin{align*}
%\lim_{n\to +\infty} \frac{s_n}{n+\sqrt{\lambda_{n+1}}}=+\infty,
%\end{align*}
%
%\begin{remark}\label{Rem1}
% It is easy to see from \eqref{estimate-v}  that if
\begin{align}\label{limit condition}
\lim_{n\to +\infty} \frac{s_n}{n+\sqrt{\lambda_{n+1}}}=+\infty,
\end{align}
%\begin{subequations}\label{limit of en}
%\begin{align}
%\lim_{n\to +\infty} e^{-\frac{s_{n-1}}{4}+C(n-1)}=&0,\label{10a}\\
%\lim_{n\to +\infty} e^{-\frac{s_{n-1}}{4}+C(n-1)+C\sqrt{\lambda_n}}=&0.
%\end{align}
%\end{subequations}
%Therefore,   if
and $
T:=\lim_{n\to\infty} t_{n}
$
exists,  then,   it holds that
\begin{align}
\lim_{t\to T^-}   \|v[t]\|_{L^2(0,1)}= 0 \quad\text{and}\quad
 \lim_{t\to T^-}V(t)= 0,\label{12}
\end{align}
which give the fixed-time stability of system~\eqref{sub-OS-1} in the fixed time  $T$.

%\end{remark}
\end{lemma}

%\begin{remark}

 It is worth noting that, under the transformation
\begin{align*}
  \widetilde{v}(x,t):=v(x,t)-\int_0^xk_n (x,y)v(y,t)\diff{y},
 \end{align*}
 system~\eqref{sub-OS-1} is converted into (see \cite{Coron:2017})
\begin{subequations}\label{target system-1}
\begin{align}
	 \widetilde{v}_{t}(x,t)=& (a(x)\widetilde{v}_{x}(x,t))_x-\lambda_n\widetilde{v}(x,t), 
   (x,t)\in [0,1]\times (t_n,t_{n+1}),\label{sub-OS-1a'}\\
	\widetilde{v}(0,t)=&0, t\in (t_n,t_{n+1}),\label{sub-OS-1b'}\\
	 \widetilde{v} (1,t)=&0  ,t\in (t_n,t_{n+1}), \label{sub-OS-1c'}
	  %\widetilde{v}(x,t_n):=&v_{0n}(x),x\in(0,1),
\end{align}
\end{subequations}
with
$
  \widetilde{v}(x,t_0):=u_0(x)-\int_0^xk_n (x,y)u_0(y)\diff{y}.
$
Since system~\eqref{target system-1} admits a unique  (weak) solution, which is defined  in the sense of \cite[(4.3) and (4.4)]{Coron:2017} and belongs to $C([t_n,t_{n+1});L^2(0,1))\cap C^1((t_n,t_{n+1});L^2(0,1))$, by virtue of Lemma~\ref{k-l}, system~\eqref{sub-OS-1} also admits a unique  (weak) solution belonging to  $C([t_n,t_{n+1});L^2(0,1))\cap C^1((t_n,t_{n+1});L^2(0,1))$.

%\end{remark}

\subsection{Control design in Case I}\label{Sec: FTS}
%In this section, we design boundary control to stabilize system~\eqref{original system} in fixed-time according to the values of the Riemann zeta function.

First,  recall that  the Riemann zeta function $\zeta(s)$  is given by (in the field of real numbers)
\begin{align*}
  \zeta(s):=\sum_{i=1}^{\infty}\frac{1}{i^{s}} ,\forall s>1,
\end{align*}
which converges for any fixed $s>1$.

  By virtue of Lemma~\ref{lem2}  and the definition of Riemann zeta function, we set
\begin{subequations}\label{choice of tn}
\begin{align}
p>&1,\\
t_0:=&0,\\
t_{n+1}-{t_n}:=& \frac{1}{(n+1)^p},n\in \mathbb{N},\label{choice of tn-b}\\
\lambda_n:=&n^{2(p+1)} +\lambda_0,n\in \mathbb{N}.
\end{align}
 \end{subequations}
It is clear that
\begin{align*}
 \frac{s_n}{n+\sqrt{\lambda_{n+1}}}\geq&  \frac{\lambda_{n-1}(t_{n}-{t_{n-1}})}{n+\sqrt{\lambda_{n+1}}}\notag\\
 \geq & \frac{\frac{(n-1)^{2(p+1)}}{n^p}}{n+\sqrt{2}(n+1)^{p+1}}\notag\\
 \geq & \frac{(n-1)^{2(p+1)}}{2\sqrt{2}(n+1)^{2p+1}}
 ,\forall n\geq \max\{1,\lambda_0\},
\end{align*}
which implies that the condition~\eqref{limit condition} is fulfilled.

For $p>1$,  it holds that
%\begin{align*}
%  \sum_{i=0}^{\infty}\frac{1}{(i+1)^{p}}=\zeta(p).
%\end{align*}
%Thus,
\begin{align*}
T_0:= t_0+\sum_{i=0}^{\infty}(t_{i+1}-{t_i})= \sum_{i=0}^{\infty}\frac{1}{(i+1)^p}=\zeta(p).
\end{align*}
%is   determined by $\zeta(p)$.

 For any constant $\sigma>0$, let $k$ be the unique solution to the kernel equation~\eqref{kernel-k} with \begin{align*}
 \lambda:=\sigma.
  \end{align*}

 Consider the disturbed system with zero initial datum
\begin{subequations}\label{sub-OS-2}
\begin{align}
	 w_{t}(x,t)=& (a(x)w_{x}(x,t))_x+ c(x)w(x,t) +f
  (x,t), (x,t)\in Q_{T_0},\label{sub-OS-2a}\\
	 w(0,t)=&d_0(t), t\in (0,T_0),\label{sub-OS-2b}\\
	 w (1,t)=&W(t) +d_1(t)  ,t\in (0,T_0), \label{sub-OS-2c}\\
	 w(x,0)=&0,x\in(0,1),
\end{align}
\end{subequations}
where
\begin{align}\label{control law-W}
  W(t):= \int_0^1k (1,y)w(y,t)\diff{y},\forall t\in (0,T_0).
 \end{align}

 To stabilize   system~\eqref{original system} in the  fixed time $T_0$,   we define the control law:
 \begin{align}\label{control law}
 U(t):=V(t)+W(t),\forall t\in (t_n,t_{n+1}),\forall n\in  \mathbb{N},
 \end{align}
where $V(t)$ and $W(t)$ are determined by \eqref{control law-V} and \eqref{control law-W}, respectively.

We state the following theorem, which is the first main result obtained in this paper, and its proof is provided in Section~\ref{Sec. V}.

\begin{theorem}\label{main result-1} Let $p$, $\{t_n\}$, and $\{\lambda_n\}$ be determined by \eqref{choice of tn}.
 Under the control law \eqref{control law},  system~\eqref{original system} admits the following   properties:
 \begin{enumerate}[(i)]%system~\eqref{original system}
 \item  the disturbance-free system~\eqref{original system}, i.e., system~\eqref{original system} with $f\equiv d_0\equiv d_1\equiv0$, is fixed-time  stable in  the time $T_0  :=\zeta(p)$, i.e.,
\begin{align*}
\lim_{t\to T_0^-}   \|u[t]\|_{L^2(0,1)}=&0;
%\lim_{t\to T_-}V(t)=&0.
\end{align*}
\item the disturbed system~\eqref{original system}, i.e., system~\eqref{original system} with $f^2+d_0^2+d_1^2\not\equiv0$, is  ISS, having the following estimate:
 \begin{align*}
  \|u[t]\|_{L^2(0,1)}
           \leq&  C_0e^{-t} \|u_0\|_{L^2(0,1)}+C_1\left( \|f\|_{L^{\infty}{(Q_{t})}}   +\|d_0\|_{{L^{\infty}(0,t)}} +\|d_1\|_{{L^{\infty}(0,t)}}\right),  \forall t\in (0,T_0),
   \end{align*}
  where $C_0$ and $C_1$ are positive constants independent of $ n$,   $u_0$, $f$, $d_0$, and $d_1$.
\end{enumerate}
\end{theorem}
%\begin{remark}

 It is worth noting that, under the transformation
\begin{align*}%\label{transformation-k}
  \widetilde{w}(x,t):=w(x,t)-\int_0^xk (x,y)w(y,t)\diff{y},
 \end{align*}
 system~\eqref{sub-OS-2} is converted into
\begin{subequations}\label{target system-2}
\begin{align}
	 \widetilde{w}_{t}(x,t)=& (a(x)\widetilde{w}_{x}(x,t))_x-\sigma  \widetilde{w}(x,t) +\widetilde{f}(x,t) , (x,t)\in Q_{T_0},\label{TS2-a}\\
	 \widetilde{w}(0,t)=&d_0(t), t\in (0,T_0),\label{TS2-b}\\
	 \widetilde{w} (1,t)=&d_1(t)   ,t\in (0,T_0), \label{TS2-c}\\
	 \widetilde{w}(x,0)=&0,x\in(0,1),
\end{align}
\end{subequations}
where
%\begin{align*}
 $\widetilde{f}(x,t) := f(x,t)-\int_0^xk (x,y)f(y,t)\diff{y}.$
  %d_0(t) :=&d_0(t) ,\\
%   d_1(t) :=&d_1(t)\left(1-\int_0^1k (1,y) \diff{y}\right).
 %\end{align*}
Since system~\eqref{target system-2} admits a unique  (weak) solution, which is defined  in the sense of \cite[(4.3) and (4.4)]{Coron:2017} and belongs to $C([0,T_0);L^2(0,1))\cap C^1((0,T_0);L^2(0,1))$, by virtue of Lemma~\ref{k-l},  system~\eqref{sub-OS-2} also admits a unique  (weak) solution belonging to  $C([0,T_0);L^2(0,1))\cap C^1((0,T_0);L^2(0,1))$. Therefore, the disturbed system~\eqref{original system} in closed loop admits a unique  (weak) solution belonging to  $C([0,T_0);L^2(0,1))$ with  $u_t[t]\in  L^2(0,1) $ for any fixed $t\in (t_n,t_{n+1})$, $n\in \mathbb{N}$.

 It is also worth mentioning that  the regularity (in the spatial variable) of {the} solution to system~\eqref{original system} can not be improved due to the fact that the initial datum of system~\eqref{target system-1} (or  system~\eqref{sub-OS-1}) is only in $L^2$-space.
However,  the regularity of   {the} solution  to system~\eqref{sub-OS-2} (or system~\eqref{target system-2})   can be improved due to the fact that   its initial datum and external disturbances
 are all continuous. For instance,  by taking spatial $H^1$-approximations of external disturbances,  it can be shown that the solutions to system~\eqref{target system-2} and system~\eqref{sub-OS-2}  belong  to $C([0,T_0);L^{\infty}(0,1))$.
% \end{remark}

\subsection{Control design in Case II}

In this section, we design a boundary controller  to stabilize system~\eqref{original system} in arbitrarily given time. For this, we let $T_0>0$ be freely prescribed and independent of initial
conditions. By virtue of Lemma~\ref{lem2},   we set
 \begin{subequations}\label{choice of tn-2}
\begin{align}
%t_0:=&0,\\
 {t_n}:=&T_0-\frac{T_0}{ n+1  },n\in \mathbb{N},\\
\lambda_n:= & n ^6 +\lambda_0,n\in \mathbb{N}.
\end{align}
  \end{subequations}

Consider \eqref{sub-OS-1} and \eqref{sub-OS-2} with $V(t)$ and $W(t)$ being determined by \eqref{control law-V} and \eqref{control law-W}, respectively.  To stabilize   system~\eqref{original system} in the freely prescribed finite time $T_0$, we define the control law by \eqref{control law}.

Note that for  $ n\geq \max\{1,\lambda_0\}$ we have
\begin{align*}
 \frac{s_n}{n+\sqrt{\lambda_{n+1}}}\geq  &\frac{\lambda_{n-1}(t_{n}-{t_{n-1}})}{n+\sqrt{\lambda_{n+1}}}
 \geq  \frac{T(n-1)^6}{2\sqrt{2} (n+1)^5 }   .
\end{align*}
Thus, the condition~\eqref{limit condition} is fulfilled    and hence   \eqref{12} holds true.

Analogous to Theorem~\ref{main result-1}, we have the following  result, which can be proved in  the same  way as  Theorem~\ref{main result-1}, and its proof is given in Section~\ref{Sec. V}.
 \begin{theorem}\label{main result-2}
 Let $\{t_n\}$ and $\{\lambda_n\}$ be determined by \eqref{choice of tn-2}.
 Under the control law \eqref{control law}, system~\eqref{original system} admits the following  properties:
  \begin{enumerate}[(i)]%system~\eqref{original system}
 \item

  the disturbance-free system~\eqref{original system}, i.e., system~\eqref{original system} with $f\equiv d_0\equiv d_1\equiv0$,  is fixed-time  stable in the freely prescribed finite time $T_0 $;
 \item the disturbed system~\eqref{original system}, i.e., system~\eqref{original system} with $f^2+d_0^2+d_1^2\not\equiv0$, is  ISS.
 \end{enumerate}
\end{theorem}
\begin{remark}\label{Remar2} Although Theorem~\ref{main result-1} and Theorem~\ref{main result-2} theoretically guarantee  the fixed-time stability of the disturbance-free system, as well as the ISS of the disturbed system, it is not feasible for practical problems to stabilize the system by designing  controllers over infinite  time intervals. Nevertheless,  within  an allowable error of fixed-time  decay estimate, the convergence estimates presented in \eqref{estimate-v}  can provide  a numeric computation scheme to implement the boundary controller over finite time intervals.
%In other words, for any $\varepsilon>0$, \eqref{estimate-v} ensures that there exists a positive integer $n_0$ and a boundary controller  $V(t)$ defined for $t\in (t_n,t_{n+1})$ with $n=0,1,...,n_0-1$  such that
%\begin{align*}
% \|u[t]\|_{L^2(0,1)}\leq \varepsilon \|u_0\|_{L^2(0,1)},\forall t\in (t_{n_0-1},t_{n_0}),
% \end{align*}
% where $u$ and $u_0$ denote the solution and initial datum of the closed-loop disturbance-free system, respectively.
\end{remark}
 \section{Stability assessment}\label{Sec. V}
 In order to prove Theorem~\ref{main result-1} and Theorem~\ref{main result-2}, we present a result on the $L^\infty$-estimate of the solution to system~\eqref{target system-2}, which plays a {crucial} role in establishing the ISS w.r.t. Dirichlet boundary disturbances for the original system~\eqref{original system} in closed loop.

 \begin{lemma}\label{lem3} {The solution of} system~\eqref{target system-2} admits the following estimate:
 \begin{align}
    \|\widetilde{w}[t]\|_{L^\infty (0,1)}
   \leq & \frac{1}{\sigma}\|\widetilde{f}\|_{L^{\infty}({Q_{t})}}+\|d_0\|_{{L^{\infty}(0,t)}} +\|d_1\|_{{L^{\infty}(0,t)}} ,  \forall t\in (0,T_0).\label{maximum estimate}
\end{align}

\end{lemma}

\begin{pf*}{Proof.} To deal with   boundary terms in the Lyapunov arguments, following \cite{Zheng:2024}, we use the generalized Lyapunov method to prove the $L^\infty$-estimate of the solution. We {first} define Stampacchia's truncation functions:
%\begin{subequations}\label{gG}
\begin{align*}
  g(s ):= &\begin{cases}\ln(1+s^2), & s \in \mathbb{R}_{> 0}, \\
0, & s\in \mathbb{R}_{\leq 0},\end{cases}  \\
G(s ):=&\int_0^s  g(\tau) \mathrm{d} \tau,\forall s\in \mathbb{R}.
\end{align*}
%\end{subequations}
It is easy to see that
\begin{subequations}\label{properties}
\begin{align}
0\leq &g(s )\leq 2|s|  ,     \forall s  \in \mathbb{R},\label{gG1} \\
g'(s)\geq &0, G(s ) \geq 0,  \forall s  \in \mathbb{R},\label{gG2} \\
g(s )  =&G(s )=0, \forall s \in \mathbb{R}_{\leq 0}.\label{gG3}
%G(\theta ) & =\frac{1}{r+1} g(\theta ) \theta , \forall\theta  \in \mathbb{R} .\label{gG3}
\end{align}
\end{subequations}
{For any $T\in (0,T_0)$, let} \begin{align*}
\Omega :=\max\left\{\frac{1}{\sigma}\|\widetilde{f}\|_{L^{\infty}({Q_{T})}},\|d_0\|_{L^{\infty}(0,T)},\|d_1\|_{L^{\infty}(0,T)} \right\}.
\end{align*}

In the generalized Lyapunov arguments, we shall
derive {$\frac{\text{d}}{\text{d} t} \int_0^1 G\left(\widetilde{w}-\Omega\right) \text{d} x$ and $\frac{\text{d}}{\text{d} t} \int_0^1 G\left(-\widetilde{w}-\Omega\right) \text{d} x$}, respectively.

 Indeed,   the boundary and initial conditions, the definition of $\Omega$, and  \eqref{gG3} imply that
\begin{align}\label{BCs}
g\left(\widetilde{w}(0, t)-\Omega\right)=g\left(\widetilde{w}(1, t)-\Omega\right)=0,\forall t\in {[0,T)}.
\end{align}
Note   that the following implication holds true for all $t\in (0,T)$:
\begin{align*}
\widetilde{w}(\cdot,t)-\Omega\geq 0 \Rightarrow    {\sigma}  \widetilde{w}(\cdot,t)\geq \|\widetilde{f}\|_{L^{\infty}(Q_{T})} \geq \widetilde{f}(\cdot,t).
\end{align*}
Therefore, it holds that
\begin{align}
\int_0^1 g\left(\widetilde{w}-\Omega\right)\left( \widetilde{f}-{\sigma}  \widetilde{w} \right) \text{d} x
   \leq  0, \forall t \in (0,T).\label{negativity}
\end{align}
Then, we deduce by \eqref{gG2}, \eqref{BCs},  and \eqref{negativity} that
\begin{align*}
 \frac{\text{d}}{\text{d} t} \int_0^1 G\left(\widetilde{w}-\Omega\right) \text{d} x 
 =&\int_0^1 g\left(\widetilde{w}-\Omega\right) \widetilde{w}_t \text{d} x \notag\\
 =&\int_0^1 g\left(\widetilde{w}-\Omega\right)\left(\left(a\widetilde{w}_{x}\right)_x -\sigma \widetilde{w} + \widetilde{f} \right)\text{d} x\notag\\
  =& -\int_0^1 ag^{\prime}\left(\widetilde{w}-\Omega\right) \widetilde{w}_x^2  \text{d} x  +\int_0^1 g\left(\widetilde{w}-\Omega\right)\left( \widetilde{f}-{\sigma}  \widetilde{w} \right) \text{d} x \notag\\
   \leq &0, \forall t \in (0,T),
\end{align*}
which implies that
\begin{align*}
 \int_0^1 G\left(\widetilde{w}(x,t)-\Omega\right) \text{d} x\leq   \int_0^1 G\left(\widetilde{w}(x,0)-\Omega\right) \text{d} x 
   \leq  0, \forall t \in (0,T).
\end{align*}
 Therefore,  $G\left(\widetilde{w} -\Omega\right)\equiv 0$ in $ Q_{T}$ and thus, for every $t\in (0,T)$, we have
  \begin{align*}
       \widetilde{w}(\cdot,t)   \leq \Omega\ \ \text{a.e.\ in}\  (0,1).
\end{align*}
  Applying similar arguments, we obtain
  \begin{align*}
 \frac{\text{d}}{\text{d} t} \int_0^1 G\left(-\widetilde{w}-\Omega\right) \text{d} x
    \leq 0, \forall t \in (0,T).
\end{align*}
  Furthermore,
  for every $t\in (0,T)$, we have
   \begin{align*}
     -\widetilde{w}(\cdot,t)   \leq \Omega\ \ \text{a.e.\ in}\  (0,1).
\end{align*}
 Finally,  for every $t\in (0,T)$, we conclude  that
    \begin{align*}
      |\widetilde{w}(\cdot,t) | \leq \Omega\ \ \text{a.e.\ in}\  (0,1),
\end{align*}
which {along with the arbitrariness of $t$ and $T$} implies \eqref{maximum estimate}.
$\hfill\blacksquare$
\end{pf*}

\begin{pf*}{Proof of Theorem~\ref{main result-1} and Theorem~\ref{main result-2}.}
The statements of  Theorem~\ref{main result-1}(i) and {Theorem~\ref{main result-2}(i)} are guaranteed by Lemma~\ref{lem2}. Therefore, it suffices to   prove the statements of  Theorem~\ref{main result-1}(ii) and Theorem~\ref{main result-2}(ii), which can be proceeded in a unified way.

Let $T_0$ be determined by the Riemann zeta function   in  Theorem~\ref{main result-1}, or freely prescribed  in Theorem~\ref{main result-2}.
In the sequel, denoted by $C$ (or $C'$)  a positive constant, which may be different from each other when appearing in different places.
 In view of Lemma~\ref{lem3} and Lemma~\ref{estimate-k}, we have
 \begin{align}
    \|\widetilde{w}[t]\|_{L^2 (0,1)}
   \leq &
    \|\widetilde{w}[t]\|_{L^\infty (0,1)}\notag\\
   \leq & \frac{1}{\sigma}\|\widetilde{f}\|_{L^{\infty}(Q_{t})}+\|d_0\|_{L^{\infty}(0,t)} +\|d_1\|_{L^{\infty}(0,t)}\notag\\
   \leq &  \frac{1}{\sigma}\left(1+e^{C\sqrt{\sigma}}\right)\|f\|_{L^{\infty}(Q_{t})} +\|d_0\|_{L^{\infty}(0,t)} +\|d_1\|_{L^{\infty}(0,t)} ,  \forall t\in (0,T_0). \label{2-norm}
\end{align}
By \eqref{2-norm}, and applying Lemma~\ref{estimate-l} and Lemma~\ref{k-l} with $\lambda:=\sigma$, we deduce that
 \begin{align}
    \|w[t]\|_{L^2 (0,1)}\leq  &
     \left(1+\|l\|_{L^{\infty}(D)}\right)\|\widetilde{w}[t]\|_{L^2 (0,1)}\notag\\
    \leq &  \left(1+\|l\|_{H^1(D)}\right)\|\widetilde{w}[t]\|_{L^2 (0,1)}\notag\\
      \leq &  \left(1+C\sigma^2\right)\|\widetilde{w}[t]\|_{L^2 (0,1)}\notag\\
        \leq &C'\left( \|f\|_{L^{\infty}(Q_{t})}+\|d_0\|_{L^{\infty}(0,t)}  +\|d_1\|_{L^{\infty}(0,t)}\right) ,  \forall t\in (0,T_0),  \label{w-2-norm}
\end{align}
where $C$ and $C'$ are positive constants independent of   $ n$,   $u_0$, $f$, $d_0$, and $d_1$.

 We infer from $u=v+w$, \eqref{w-2-norm}, and Lemma~\ref{lem2} that
\begin{align}
 \|u[t]\|_{L^2(0,1)}
 \leq& \|v[t]\|_{L^2(0,1)}+\|w[t]\|_{L^2(0,1)}\notag\\
 \leq& Ce^{-\frac{s_{n-1}}{4}+C(n-1)} \|u_0\|_{L^2(0,1)} +C'\left( \|f\|_{L^{\infty}(Q_{t})}+\|d_0\|_{L^{\infty}(0,t)}
  +\|d_1\|_{L^{\infty}(0,t)}\right) ,\forall t\in (t_n,t_{n+1}),\label{Equ.20}
    %  \leq&  2(1+C \sigma^2)\frac{(1+e^{ C\sqrt{\sigma}})}{ \sqrt{\varepsilon(2\sigma -\varepsilon)}} \sup_{s\in (0,T_0)} \|f(\cdot,s)\|_{L^2(0,1)}\notag\\
% &+ Ce^{-\frac{s_{n-1}}{4}+C(n-1)} \|u_0\|_{L^2(0,1)}^2\notag\\
% &+Ce^{-\frac{s_{n-1}}{4}+C(n-1)},\forall t\in (t_n,t_{n+1}),
   \end{align}
    where $\{s_n\}$ is defined by \eqref{sn},   and $C,C'$ are positive constants independent of $ n$,   $u_0$, $f$, $d_0$, and $d_1$.

 In view of \eqref{limit condition}, there exists a sufficiently large  integer $N_0\in \mathbb{N}_{+}$ such that
    \begin{align}\label{T-condition}
  \frac{s_{n-1}}{4}-C(n-1) >T_0,\forall n>N_0.
   \end{align}
   Moreover, for any $t\in (0,T_0)$,  there exists a positive integer $n_0$ such that $t\in (t_{n_0},t_{n_0+1})$ with either $n_0>N_0$ or $n_0+1\leq N_0$.

   If $n_0>N_0$, we deduce by \eqref{Equ.20}  and \eqref{T-condition}  that
            \begin{align}
 \|u[t]\|_{L^2(0,1)}
       \leq&  Ce^{-T} \|u_0\|_{L^2(0,1)}  +C'\left( \|f\|_{L^{\infty}(Q_{t})} +\|d_0\|_{L^{\infty}(0,t)}  +\|d_1\|_{L^{\infty}(0,t)}\right)\notag\\
 \leq&
 Ce^{-t} \|u_0\|_{L^2(0,1)} +C'\left( \|f\|_{L^{\infty}(Q_{t})} 
       +\|d_0\|_{L^{\infty}(0,t)} +\|d_1\|_{L^{\infty}(0,t)}\right),\label{23}
   \end{align}
where $C$ and $C'$ are positive constants independent of   $ n$,   $u_0$, $f$, $d_0$, and $d_1$.

  If $n_0+1\leq N_0$, we infer  from  \eqref{Equ.20} that there exists a positive constant $C''$ depending only on $N_0$ such that
     \begin{align}
  \|u[t]\|_{L^2(0,1)}
        \leq&   C'' \|u_0\|_{L^2(0,1)} +C'\left( \|f\|_{L^{\infty}(Q_{t})} 
        +\|d_0\|_{L^{\infty}(0,t)} +\|d_1\|_{L^{\infty}(0,t)}\right)\notag\\
  \leq&   C''e^{t_{N_0+1}}e^{-t} \|u_0\|_{L^2(0,1)}   +C'\left( \|f\|_{L^{\infty}(Q_{t})}+\|d_0\|_{L^{\infty}(0,t)} 
       +\|d_1\|_{L^{\infty}(0,t)}\right),\label{24}
   \end{align}
where  $C'$ is positive constant  independent of   $ n$,   $u_0$, $f$, $d_0$, and $d_1$.

    Combining \eqref{23} and \eqref{24}, we conclude that disturbed system in closed loop is  ISS, having the estimate for all $t\in (0,T_0)$:
    \begin{align*}
 \|u[t]\|_{L^2(0,1)}
           \leq&
C \left(e^{-t} \|u_0\|_{L^2(0,1)}+  \|f\|_{L^{\infty}(Q_{t})} 
 +\|d_0\|_{L^{\infty}(0,t)} +\|d_1\|_{L^{\infty}(0,t)}\right),
   \end{align*}
   where $C$ is a positive constant  independent of   $ n$,   $u_0$, $f$, $d_0$, and $d_1$.
$\hfill\blacksquare$
\end{pf*}
%\begin{remark}
% It is worth noting that, for  any $T_0>0$, system \eqref{original system} is input-to-state stable under the common boundary control law
% \begin{align}\label{control law-ISS}
%  U(t):= \int_0^1k (1,y)u(y,t)\diff{y},\forall t\in (0,T_0),
% \end{align}
% where $k$ is determined by \eqref{kernel-k}. Nevertheless, system \eqref{original system} is not fixed-time stable when $f\equiv 0$.
%\end{remark}

\section{Numerical example}\label{Sec. VI}
We consider the following {system}
\begin{subequations}\label{original system-simulation}
\begin{align}
	 u_{t}(x,t)=&  u_{xx}(x,t) + cu(x,t) ,(x,t)\in Q_{T_0}, \\
	 u(0,t)=&0, t\in (0,T_0), \\
	 u (1,t)=&U(t)+d_1(t)   ,t\in (0,T_0),
\end{align}
\end{subequations}
with
  %\begin{subequations}
%\begin{align*}
$c= 24$, $
d_1(t)=  A\sin(30t)$, and $A\in\{0,1,2\}$,
%\end{align*}
 %\end{subequations}
 where $A$ is used to describe {the amplitude of different}  Dirichlet  boundary disturbances.

 The initial condition is given by
 \begin{align*}
u(x,0)=&u_0(x),x\in(0,1),
\end{align*}
or,
  \begin{align*}
 u (x,0)=&10u_0(x),x\in(0,1),
\end{align*}
where
 %\begin{align*}
$u_0(x):=   -4 \sin(15(x-0.5 ) )$.
%\end{align*}

\begin{figure}[htbp]
\begin{center}
%\subfigure{\includegraphics[scale=0.65]{opensol.eps}\label{Fig.opensol}}\\
%(a) Trajectory of the solution\\
 \includegraphics[scale=0.65]{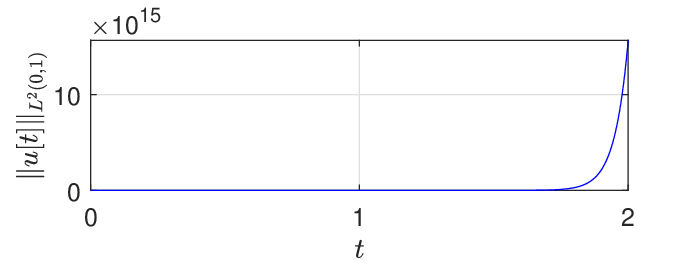}\label{Fig.opennorm}
%$L^2$-norm of the solution\\
\caption{$L^2$-norm of  solution  to the open-loop   system~\eqref{original system-simulation}   with $d_1\equiv 0$  and initial datum $u_0(x)$} \label{Fig.openloop}
\end{center}
\end{figure}

 It is known that the  disturbance-free system~\eqref{original system-simulation} in open loop  is unstable due to the fact that $c>\pi^2$; see also  Fig.~\ref{Fig.openloop}

% It is worth noting that, for  any $T_0>0$, system \eqref{original system} is input-to-state stable under the common boundary control law
% \begin{align}\label{control law-ISS}
%  U(t):= \int_0^1k (1,y)u(y,t)\diff{y},\forall t\in (0,T_0),
% \end{align}
% where $k$ is determined by \eqref{kernel-k}. Nevertheless, system \eqref{original system} is not fixed-time stable when $f\equiv 0$.

  In order to input-to-state stabilize system~\eqref{original system-simulation}, we split it into two subsystems:
%\begin{subequations}
\begin{align*}
	 v_{t}(x,t)=&  v_{xx}(x,t) + c v(x,t), (x,t)\in (0,1)\times (t_n,t_{n+1}), \\
	v(0,t)=&0, t\in (t_n,t_{n+1}), \\
	 v (1,t)=&V(t)   ,t\in (t_n,t_{n+1}), \\
	 % v(x,t_n):=&v_{0n}(x),x\in(0,1),\\
v(x,t_0) =&u(x,0), x\in(0,1),
\end{align*}
%\end{subequations}
and
%\begin{subequations}
\begin{align*}
	 w_{t}(x,t)=&  w_{xx}(x,t) + c w(x,t) , (x,t)\in Q_{T_0}, \\
	 w(0,t)=&0, t\in (0,T_0), \\
	 w (1,t)=&W(t)+d_1(t)   ,t\in (0,T_0), \\
	 w(x,0)=&0,x\in(0,1),
\end{align*}
%\end{subequations}
where
\begin{align*}
V (t)=& \int_0^1k_n (1,y)v(y,t)\diff{y},\forall t\in (t_n,t_{n+1}),\\
  W(t) =& \int_0^1k (1,y)w(y,t)\diff{y},\forall t\in (0,T_0).
   \end{align*}
   According to Remark~\ref{Remark1}, the kernel functions are given by
   \begin{align*}
   k_n(1,y)=&- (\lambda_n+c) y\frac{I_1\left(\sqrt{(\lambda_n+c)(1-y^2)}\right)}{\sqrt{(\lambda_n+c)(1-y^2)}},\\
   % \end{align*}
%   with $\lambda_n\geq \lambda_0$, and
%      \begin{align*}
   k (1,y)=&- (\sigma+c) y\frac{I_1\left(\sqrt{(\sigma+c)(1-y^2)}\right)}{\sqrt{(\sigma+c)(1-y^2)}}.
 \end{align*}

According to Theorem~\ref{main result-1} and Theorem~\ref{main result-2},  the following control law
 \begin{align}\label{final control}
  U(t):=V(t)+W(t),\forall t\in (t_n,t_{n+1})
 \end{align}
 can be used to input-to-state stabilize system~\eqref{original system-simulation}
after choosing appropriate $t_n, \lambda_n\geq \lambda_0$, and $\sigma$.

In simulations, we always choose
\begin{align*}
\lambda_0:= 3.5 \quad \text{and}\quad
 \sigma=   {1}.
\end{align*}
 \textbf{Numerical results in Case I.} In this case, we choose $p=1.9$, then the fixed time is given by
 \begin{align*}
 T_0= \zeta(p)=1.7497,
 \end{align*} and
 $t_n,\lambda_n$ are given by
 \begin{align*}
 t_0:=&0,\\
t_{n+1}  =& {t_n}+\frac{1}{(n+1)^p},{n\in \mathbb{N}_+},\\
\lambda_n =& n^{2(p+1)} +\lambda_0, n\in \mathbb{N} .
 \end{align*}
 By virtue of Remark~\ref{Remar2}, for   numeric computations, we   consider {$n=0,1,2$}. Thus, the time interval $(0,T_0)$ is divided into
 \begin{align*}
(t_0,t_1)=&(0,1),\\
(t_1,t_2)=&(1, 1.2679),\\
(t_2,t_3)=&(1.2679, 1.7497),
\end{align*}
for which $\{\lambda_n\}$ is given by
 \begin{align*}
 {\lambda_0=}    {3.5}, \quad
 \lambda_1=   4.5, \quad \text{and}  \quad
\lambda_2=   58.2152,
\end{align*}
 respectively.

\begin{figure}[htbp]
\begin{center}
\subfigure{\includegraphics[scale=0.65]{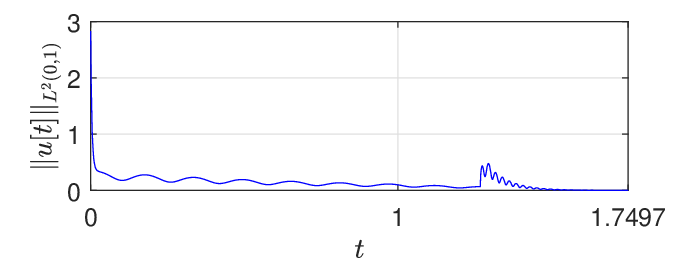}}\\
(a) $u(x,0)=u_0(x),A=0$\\
\subfigure{\includegraphics[scale=0.65]{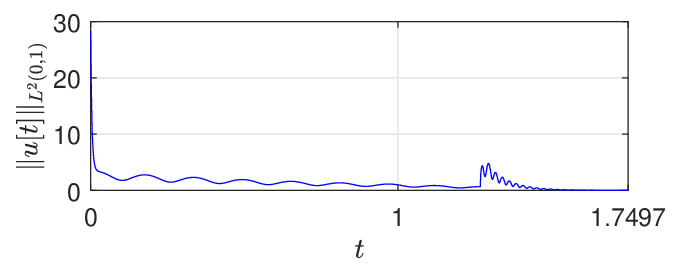}}\\
(b) $u(x,0)=10u_0(x),A=0$\\
\caption{{$L^2$-norms of the solutions to the disturbance-free closed-loop system with different initial data} under the control law \eqref{final control} in Case I}\label{Fig.FTS1}
\end{center}
\end{figure}

  {Figure~\ref{Fig.FTS1}~(a)} shows that system~\eqref{original system-simulation} with   initial datum $u_0$ is fixed-time stable in $L^2$-norm in the time $T_0=1.7497$,  while   {Fig.~\ref{Fig.FTS1}~(b)} shows that the system is still fixed-time stable in the same time $T_0$ when the initial datum is increased to be $10u_0$.

\begin{figure}[htbp]
\begin{center}
\subfigure{\includegraphics[scale=0.65]{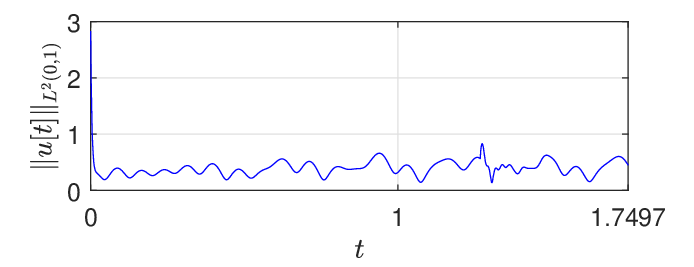}}\\
(a) $u(x,0)=u_0(x),A=1$\\
\subfigure{\includegraphics[scale=0.65]{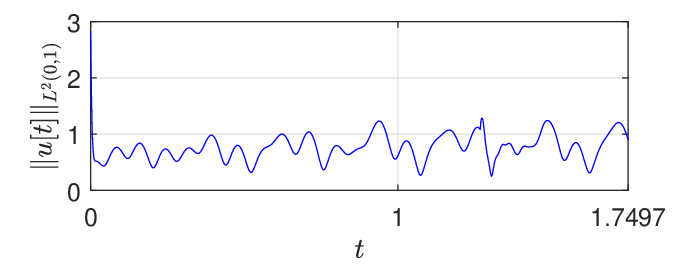}}\\
(b)  $u(x,0)=u_0(x),A=2$\\
\caption{{$L^2$-norms of  the solutions}    to the closed-loop system \eqref{original system-simulation} with different  disturbances  under the control law \eqref{final control} in Case I}\label{Fig.FTISS1}
\end{center}
\end{figure}

   Under the same initial condition,  {Fig.~\ref{Fig.FTS1}~(a) and  Fig.~\ref{Fig.FTISS1}} show that the $L^2$-norms of {the} solutions to system~\eqref{original system-simulation} with different boundary disturbances remains bounded due to the facts that the  disturbances are bounded. However, the  norms are  more influenced by disturbances with large  amplitudes. This reflects well the ISS property of the system in the closed loop.

 %\begin{subequations}
%\begin{align}
%p=&1.9,\\
%t_0:=&0,\\
%t_{n+1} :=&{t_n}+\frac{1}{(n+1)^p},\forall n\geq 0,
%\\
% T_0=&\zeta(p)=1.7497,\\
%(t_0,t_1)=&(0,1),\\
%(t_1,t_2)=&(1, 1.2679),\\
%(t_1,t_2)=&(1.2679, 1.7497),\\
%\end{align}
% \end{subequations}
%  \begin{subequations}
%\begin{align}
%{\lambda_0:=}&  {3.5},\\
%\lambda_n:=&n^{2(p+1)} +\lambda_0,\forall n\geq 0,\\
%n=&2,\\
%\lambda_1:=&  4.5,\\
%\lambda_2:=&  58.2152.
%\end{align}
% \end{subequations}

 %\begin{figure}[htbp]
%\begin{center}
%\subfigure{\includegraphics[scale=0.65]{FTISSA0sol.eps}}\\
%(a) Trajectory of the solution when $A=0$\\
%\subfigure{\includegraphics[scale=0.65]{FTISSA1sol.eps}}\\
%(b) Trajectory of the solution when $A=1$\\
%\subfigure{\includegraphics[scale=0.65]{FTISSA2sol.eps}}\\
%(c) Trajectory of the solution when $A=2$\\
%\caption{Trajectories of  solutions  to the closed-loop system \eqref{original system-simulation}   under the control law \eqref{final control} with $\lambda=\lambda_0=1$} \label{Fig.FTISS0}
%\end{center}
%\end{figure}


   \textbf{Numerical results in Case II.} In this case, we freely prescribe the finite time to be
\begin{align*}
T_0 =1.5.
\end{align*}
Then, $t_n,\lambda_n$ are given by
 \begin{align*}
t_0 =&0,\\
 {t_n} = &T_0-\frac{T_0}{ n+1  }, n\in \mathbb{N}_+, \\
 \lambda_n =&  n ^6 +\lambda_0,n\in \mathbb{N}.
 \end{align*}
 By virtue of Remark~\ref{Remar2},  for   numeric computations, we   consider {$n=0,1,2$}. Thus, the time interval $(0,T_0)$ is divided into
 \begin{align*}
 (t_0,t_1)=& (0,0.75),\\
 (t_1,t_2)=& (0.75,1),\\
 (t_2,t_3)=& (1,1.5),
\end{align*}
associated with
 \begin{align*}
\lambda_0 =3.5,\quad
\lambda_1=4.5,\quad \text{and}  \quad
\lambda_2= 67.5,
\end{align*}
 respectively.

\begin{figure}[htbp]
\begin{center}
\subfigure{\includegraphics[scale=0.65]{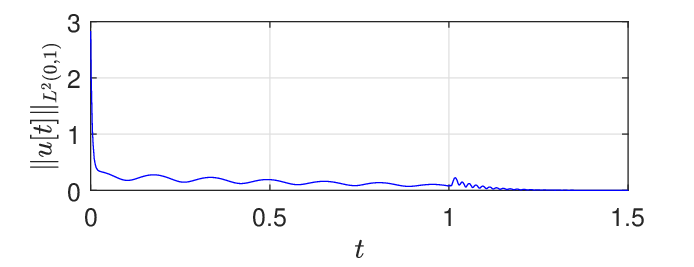}}\\
(a)  $u(x,0)=u_0(x),A=0$\\
\subfigure{\includegraphics[scale=0.65]{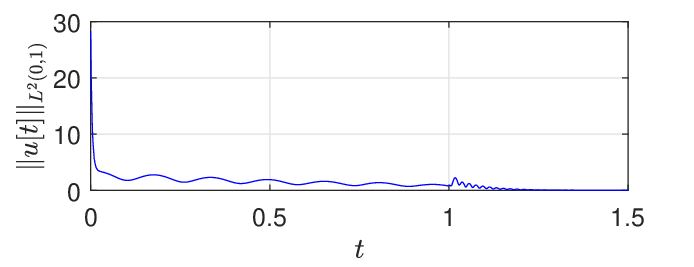}}\\
(b)$u(x,0)=10u_0(x),A=0$\\
\caption{{$L^2$-norms of  the solutions   to the disturbance-free closed-loop system} \eqref{original system-simulation} with different  initial data  under the control law \eqref{final control} in Case II}\label{Fig.PTS1}
\end{center}
\end{figure}

\begin{figure}[htbp]
\begin{center}
\subfigure{\includegraphics[scale=0.65]{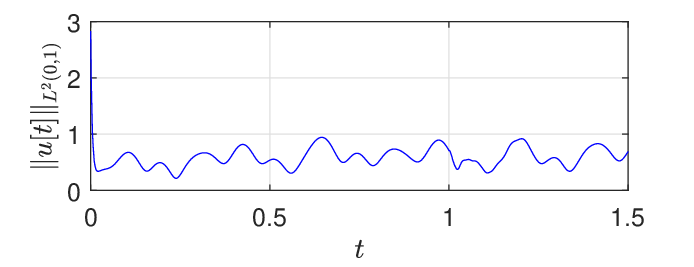}}\\
(a)  $u(x,0)=u_0(x),A=1$\\
\subfigure{\includegraphics[scale=0.65]{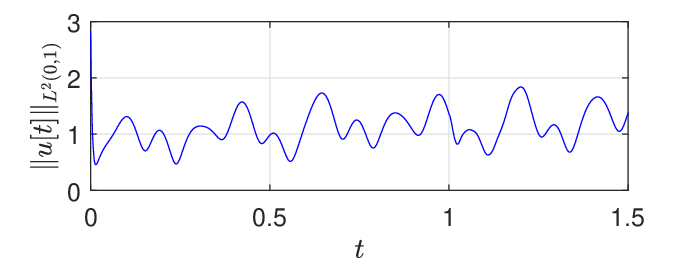}}\\
(b)  $u(x,0)=u_0(x),A=2$\\
\caption{{$L^2$-norms of  the solutions}   to the closed-loop system \eqref{original system-simulation} with different  disturbances  under the control law \eqref{final control} in Case II}\label{Fig.PTISS1}
\end{center}
\end{figure}

  {Figure~\ref{Fig.PTS1}~(a)} shows that system~\eqref{original system-simulation} with   initial datum $u_0$ is fixed-time stable in $L^2$-norm in the time $T_0=1.5$,  while   {fig.~\ref{Fig.PTS1}~(b)} shows that the system is still fixed-time stable in the same time $T_0$ when the initial datum is increased to be $10u_0$.
  As in Fig.~\ref{Fig.FTS1}~(a)  {and Fig.~\ref{Fig.FTISS1}},  under the same initial condition,  {Fig.~\ref{Fig.PTS1}~(a) and Fig.~\ref{Fig.PTISS1}} also reflect well the ISS property of the system in the closed loop.

%\begin{align}
%t_0:=&0,\\
% {t_n}:=&T_0-\frac{T_0}{ n+1  },\forall n\geq 1,\\
% (t_0,t_1)=&(0,0.75),\\
% (t_1,t_2)=&(0.75,1),\\
% (t_2,t_3)=&(1,1.5),\\
%\lambda_n:=& n ^6 +\lambda_0,\forall n\geq 0,\\
%\lambda_0:=3.5,\\
%\lambda_1=4.5,\\
%\lambda_2= 67.5.
%\end{align}

\section{Conclusion}\label{Sec. VII}
This paper addressed the problem of input-to-state stabilization for   $1$-D parabolic   equations in {the} presence of distributed in-domain and Dirichlet boundary disturbances {while ensuring} the fixed-time stability of the corresponding disturbance-free system.
 A   boundary controller  was designed by using backstepping based on the technique of splitting, and the stability was assessed by using the generalized Lyapunov method, which can be used to deal with boundary terms easily.
   Numerical simulations were conducted, and the results confirmed the ISS properties of the considered system  with in-domain and Dirichlet {boundary}  disturbances, as well as the fixed-time stability of the associated disturbance-free system.

   It is worth mentioning that the proposed scheme is suitable for input-to-state  {stabilization of}  $1$-D linear parabolic  equations with  Robin or Neumann boundary disturbances,   as well as for fixed-time {stabilization of}  the corresponding disturbance-free system, for which the classical Lyapunov method remains valid for stability analysis.
   However, it is still  challenging to directly apply  the proposed method to fixed-time  input-to-state  {stabilization of} the system  with Dirichlet boundary disturbances due to the fact that  {the fixed-time ISS is   a much stronger property}  than the ISS. This problem will be addressed in our  future work.
%\section*{Declaration of competing interest}
%The authors declare that they have no known competing financial interests or personal relationships that could have appeared to influence the work reported in this paper.
%
%\section*{Data availability}
%No data was used for the research described in the article.
%
%\section*{Acknowledgments}
%This work is supported in part by NSFC under grant NSFC-11901482 and in part by NSERC under grant {RGPIN-2024-04709}.

%\bibliographystyle{plain}        % Include this if you use bibtex
% %\bibliography{Referencesfullname}        % and a bib file to produce the
%                                 % bibliography (preferred). The
%% Loading bibliography database
%\bibliography{References}
%%%%%%%%%%%%%%%%%%%%%%%%%%%%%%%%%%%%%%%%%%%%%%

\end{document}